# GENERALIZED SCHATUNOWSKY THEOREM IN A WEAK ARITHMETIC

HALA KING AND VICTOR PAMBUCCIAN

ABSTRACT. Schatunowsky's 1893 theorem, that 30 is the largest number all of whose totatives are primes, has been recently generalized by Kaneko and Nakai. In its generalized form, it states the finiteness of the set of all positive numbers $n$, which, for a fixed prime $p$, have the property that all of $n$'s totatives that are not divisible by any prime $\leq p$ are prime numbers. It is this generalized form that we show holds in a weak arithmetic.

## 1. INTRODUCTION

Schatunowsky [1, p. 132] proved that 30 is the largest number $n$, all of whose totatives (numbers greater than 1 but less than $n$ which are relatively prime to it) are primes.

Fixing a prime number $p$, we will call *p-good* all numbers $n$ with the property that all of $n$'s totatives that are not divisible by any prime $\leq p$ are prime numbers. In [2], Schatunowsky's result has been generalized by showing that, for any given prime $p$, the set of $p$-good numbers is finite.

Proofs of elementary number-theoretic results, when considered from a logical point of view, are taken to happen inside Peano Arithmetic (PA), which contains an axiom schema for induction, in other words, a countably infinite number of axioms. The compactness theorem for first-order logic tells us that each proof of an elementary number-theoretic result $\varphi$, which can be proved inside PA, can be proved from a subset $\Sigma$ of the axiom system of PA, consisting of finitely many axioms. In other words, that it is enough to consider finitely many instances of the induction axiom schema. However, the compactness theorem is not costructive, so it does not allow for the determination of the exact finite subsystem $\Sigma$ from which $\varphi$ can be derived.

Finding such finite subsystems can, however, be done by analyzing an actual proof of $\varphi$. While the specific cases of the induction axiom needed for the proof of $\varphi$ are not illuminating, it is more interesting to find out which particular elementary number-theroretical statements are needed for that proof.

The aim of this note is to provide just such a proof of the above-mentioned generalization of Schatunowsky's theorem from a finite set of axioms, consisting of all the axioms of Peano Arithmetic except the induction axiom schema, which is replaced by a few elementary number-theoretic statements with a clear meaning.

Our base theory will be that very weak fragment of Peano Arithmetic, referred to as $PA^-$ — whose models are the positive cones of discretely ordered rings — which consists of all axioms except the induction axiom schema.

## 2. The axioms of PA$^-$

PA$^-$ is a theory expressed in a language with $+, \cdot, 1, 0, <$ as primitive notions. Its axiom system consists of 15 axioms. We will reproduce them here from [3, pp. 16-18] for the reader's convenience, and will omit the universal quantifiers for all universal axioms.

**A 1.** $(x + y) + z = x + (y + z)$

**A 2.** $x + y = y + x$

**A 3.** $(x \cdot y) \cdot z = x \cdot (y \cdot z)$

**A 4.** $x \cdot y = y \cdot x$

**A 5.** $x \cdot (y + z) = x \cdot y + x \cdot z$

**A 6.** $x + 0 = x \wedge x \cdot 0 = 0$

**A 7.** $x \cdot 1 = x$

**A 8.** $(x < y \wedge y < z) \to x < z$

**A 9.** $\neg x < x$

**A 10.** $x < y \vee x = y \vee y < x$

**A 11.** $x < y \to x + z < y + z$

**A 12.** $(0 < z \wedge x < y) \to x \cdot z < y \cdot z$

**A 13.** $(\forall x)(\forall y)(\exists z)\, x < y \to x + z = y$

**A 14.** $0 < 1 \wedge (x > 0 \to (x > 1 \vee x = 1))$

**A 15.** $x > 0 \vee x = 0$

Let PA$^-$ be the theory axiomatized by A1-A15. This theory is very weak, in which one cannot show that among two consecutive numbers one is even (or one is odd). Here by "even" and "odd" we mean that an element is a multiple of 2 or a multiple of 2 plus 1, where 2 stands for $1 + 1$. In fact, there are models of PA$^-$ in which, for any natural number $n$, there are sequences of $n$ consecutive numbers, none of which is odd or even. This can be seen by taking a look at the positive cone of $\mathbb{Z}[X]$, which is a model of PA$^-$, and where the sequence $X+1, \ldots, X+n$ has no even element and no odd element. Nor can one show in PA$^-$ that the concepts of "irreducible" (a number having no other divisors that 1 and itself) and "prime" (a number which, whenever it divides a product it must divide one of the factors) coincide. The two notions are thus defined by

(1) $\pi_1(x)\ :\Leftrightarrow\ (\forall a)(\forall b)\, 1 < x \wedge (x = a \cdot b \to (a = 1 \vee b = 1))$
(2) $\pi_2(x)\ :\Leftrightarrow\ (\forall a)(\forall b)(\forall c)(\exists d)\, 1 < x \wedge (x \cdot c = a \cdot b \to (x \cdot d = a \vee x \cdot d = b))$

and there are models of PA$^-$ in which the two are not equivalent (see [3, pp. 21]). It is plain that

(3) $$\mathrm{PA}^- \vdash (\forall x)\, \pi_2(x) \to \pi_1(x).$$

The 'primes' of the generalized Schatunowsky theorem can thus be those numbers $x$ satisfying $\pi_1(x)$ or those satisfying $\pi_2(x)$. We will choose to interpret '$x$ is prime' as $\pi_2(x)$.

For $u \geq 1$, we denote by $\overline{u}$ the term $((\ldots((1+1)+1)+\ldots)+1)$, in which there are $u$ many 1s, and we let $\overline{0}$ be 0 itself; the terms $\overline{u}$ will be referred to as *numerals*.

We will use the following abbreviations: $x \leq y :\Leftrightarrow x < y \vee x = y$, $x|y :\Leftrightarrow (\exists z) \, xz = y$, for '$x$ divides $y$', and $\varrho(m,n) :\Leftrightarrow (\forall d) \, d|m \wedge d|n \rightarrow d = 1$ for '$m$ and $n$ are relatively prime'.

## 3. Additional axioms

Given how weak $PA^-$ is, it is no surprise that we need to state four additional axioms.

The first two state that any prime number has a *successor* prime, a *next* prime, and any prime $> 2$ has a *predecessor* prime, a *previous* prime.

**A 16.** $(\forall p)(\exists q)(\forall u) \, \pi_2(p) \rightarrow (p < q \wedge \pi_2(q) \wedge \pi_2(u) \wedge p < u \rightarrow q \leq u)$

**A 17.** $(\forall p)(\exists q)(\forall u) \, \pi_2(p) \rightarrow (q < p \wedge \pi_2(q) \wedge \pi_2(u) \wedge u < p \rightarrow u \leq q)$

We will denote the $q$ in A16 by $S(p)$ and the $q$ in A17 by $P(p)$.

Our third axiom states that, for every element $n$ greater than $\overline{4}$, there exists a largest prime whose square is less than n.

**A 18.** $(\forall n)(\exists p)(\forall q) \, \overline{4} < n \rightarrow \pi_2(p) \wedge p^2 < n \wedge (p < q \wedge \pi_2(q) \rightarrow n \leq q^2)$

The next axiom we want to state is an inequality in which three consecutive primes are involved. To simplify its statement, we introduce the following defined predicate $\sigma$, which states that $a$ and $b$ are consecutive primes:

$\sigma_u(a,b) :\Leftrightarrow \pi_2(a) \wedge \pi_2(b) \wedge a < b \wedge (a < u \wedge \pi_2(u) \rightarrow b \leq u)$

We are now ready to state the axiom stating that the square of a prime $\geq 19$ is less than twice the product of the previous two primes:

**A 19.** $17 < q \wedge \sigma_u(r,p) \wedge \sigma_u(p,q) \rightarrow q^2 < 2pr$

That A19 holds in the standard model $\mathbb{N}$ follows from [2, p. 445, (6)].

## 4. The generalized Schatunowsky theorem and its proof

The generalized Schatunowsky theorem, as proved in [2], can now be stated in a weak form, by asking for the existence, for any given prime number $p$, for an upper bound among all $p$-good numbers:

**GS$^w$ .** $(\forall p)(\exists n)(\forall m)(\exists q)(\forall u) \, \pi_2(p) \wedge n \leq m \rightarrow [\varrho(q,m) \wedge q < m \wedge 1 < q$
    $\wedge (u \leq p \wedge \pi_2(u) \rightarrow \varrho(u,q)) \wedge \neg \pi_2(q)]$.

It can also be stated in a strong form, asking, for any given prime number $p$, for a largest $p$-good number:

**GS$^s$ .** $(\forall p)(\exists n)(\forall t)(\forall v)(\forall m)(\exists q)(\forall u) \, \pi_2(p) \rightarrow [(1 < t \wedge t < n$
    $\wedge \varrho(t,n) \wedge (v \leq p \wedge \pi_2(v) \rightarrow \varrho(v,t)) \rightarrow \pi_2(t)) \wedge (n \leq m \rightarrow (\varrho(q,m)$
    $\wedge q < m \wedge 1 < q \wedge (u \leq p \wedge \pi_2(u) \rightarrow \varrho(u,q)) \wedge \neg \pi_2(q)))]$.

Let $\Sigma$ stand for the theory axiomatized by the axioms A1-A19.

**Theorem 4.1.** $\Sigma \vdash \mathbf{GS}^w$

*Proof.* Let $\mathfrak{M}$ be a model of $\Sigma$. If $\mathfrak{M} = \mathbb{N}$, then $\mathbf{GS}^w$ holds by [2, Th. 4 & 7]. Suppose $\mathfrak{M}$ is a nonstandard model of $\Sigma$. Given the presence of A18 among our axioms, there must be nonstandard primes in $\mathfrak{M}$. Let $p \geq \overline{7}$ be a prime in $\mathfrak{M}$. We will prove that no number greater than or equal to $S(S(S(p)))^2 + 1$ (we recall that $S(p)$ is the succesor prime to $p$) can be $p$-good. In other words, we will prove that $\mathbf{GS}^w$ holds with $n = S(S(S(p)))^2 + 1$. Let $m \geq n$. Let $q$ be the largest prime, known to exist by A18, for which $q^2 < m$. Then $S(S(S(p))) \leq q$. Since $p \geq \overline{7}$, we have $q \geq \overline{17}$ and $S(q) \geq \overline{19}$. By A19, we know that

$$(4) \qquad S(q)^2 < \overline{2}qP(q).$$

Notice that, if $s$ is a prime with $s^2 < v$, then $v$ and $s^2$ are co-prime unless $s|v$. For, if $d$ were a divisor of both $v$ and $s^2$, we would have $dx = v$ and $dy = s^2$. Since $s$ is prime, i. e. $\pi_2(s)$ holds, and $s$ divides $dy$, we must have $s|d$ or $s|y$. If $s|d$, then $s|v$ and we are done. If $s|y$, i. e. $st = y$ for some $t$, then $dst = s^2$, so $dt = s$, which, given that $s$ is irreducible, implies that $d = s$ or $t = s$. If $d = s$, then $s|v$ and we are done. If $t = s$, then, from $st = y$, we conclude that $y = s^2$, and from $dy = s^2$ that $d = 1$.

Given that $q^2 < m$, we conclude that $q^2$ and $m$ are either co-prime or else $q \,|\, m$. If $q^2$ and $m$ were co-prime, then $m$ has a totative — namely $q^2$, which is co-prime with all the primes $\leq p$, since $q \geq S(S(S(p)))$ — that is not a prime. The same argument can be made for $P(q)$ or $P(P(q))$ instead of $q$, and since $P(q) \geq S(S(p))$ and $P(P(q) \geq S(p)$, the same conclusion holds if $P(q)^2$ or $P(P(q))^2$ were co-prime with $m$.

On the other hand, if each of $q$, $P(q)$, and $P(P(q))$ were divisors of $m$, then, since $q$, $P(q)$, and $P(P(q))$ are primes, $qP(q)P(P(q))\,|\,m$ as well. Thus $qP(q)P(P(q)) \leq m < S(q)^2 < \overline{2}qP(q)$, the last inequality by (4). That would imply $P(P(q)) < \overline{2}$, a contradiction.

If $p \in \{\overline{2}, \overline{3}, \overline{5}\}$, then use the same argument with $n = \overline{290}$ to find that, if $m \geq n$, then the greatest prime $q$ with $q^2 < m$ is $\geq \overline{17}$ and conclude that $m$ is $\overline{5}$-good. This means that there cannot be nonstandard $m$ having the desired property, so all elements with that property are those listed in [2, Th. 7]. □

To prove the validity of $\mathbf{GS}^s$, we need two additional axioms. The first one is Proposition 30 in Book VII of Euclid's *Elements*. It states that every number greater than 1 has a prime divisor. Formally

**A 20.** $(\forall n)(\exists p)\, n > 1 \rightarrow (\pi_2(p) \wedge p|n)$

The next axiom states that, for any prime $p$, there is a greatest number $< \frac{S(p)^2}{p}$. Formally:

**A 21.** $(\forall pqu)(\exists k)\, \sigma_u(p,q) \rightarrow (kp < q^2 \wedge (k+1)p > q^2)$

Let $\Sigma'$ stand for the theory axiomatized by the axioms A1-A21.

**Theorem 4.2.** $\Sigma' \vdash \mathbf{GS}^s$

*Proof.* If $p$ in $\mathbf{GS}^s$ is $< \overline{7}$, then the last paragraph in the proof of Theorem 4.1 ensures that the only $p$-good numbers are those in the Table of [2, Th. 7]. If $p = \overline{7}$, then, by the first part of the proof of Theorem 4.1, no number $\geq \overline{290}$ can be $\overline{7}$-good, so, according to the Table of [2, Th. 7], the largest $\overline{7}$-good number is $\overline{286}$.

Suppose now $p > \overline{7}$. We want to show that **GS**$^s$ holds with $n = S(p)k_p$, where $k_p$, which exists by A21, is such that $k_p S(p) < S(S(p))^2$ and $(k_p + 1)S(p) > S(S(p))^2$. To show that $n$ is $p$-good, assume there exists a non-prime totative $t$ of $n$, which is not divisible by any prime $\leq p$. By A20, $t$ has a prime divisor $q$, which, being greater than $p$ and different from $S(p)$, which is a divisor of $n$, must be $\geq S(S(p))$. Thus $t = qu$, with $q \geq S(S(p))$. Since $t$ is not prime, $u > 1$, and so $u$ has, by A20, a prime divsior $s$ as well, which, for the same reasons invoked earlier, must also be $\geq S(S(p))$. Thus $t \geq S(S(p))^2$. On the other hand, $t < n < S(S(p))^2$, a contradiction.

Suppose $m > n$. We will show that $m$ is not $p$-good. Since as $S(p) < S(S(p)) - 1$, we have $(S(S(p))+1)S(p) < S(S(p))^2$. Since $k_p$ is the largest number with $k_p S(p) < S(S(p))^2$, we have $k_p > S(S(p))$. Since $S(p)^2 < m$, if $m$ were $p$-good, then $S(p)$ would have to divide $m$, for else $S(p)^2$ would be a totative of $m$, which is not divisible by any prime $\leq p$. Thus $m = S(p)u$, for some $u > k_p$. By the definition of $k_p$, we have $S(p)u > S(S(p))^2$, and so, if $m$ were to be $p$-good, $S(S(p))$ would have to divide $m$, or else $S(S(p))^2$ would be a totative of $m$, which is not divisible by any prime $\leq p$. Since both $S(p)$ and $S(S(p))$ are prime, we conclude that their product must divide $m$, i.e., $m = S(p)S(S(p))v$. Since $n = k_p S(p) > S(S(p))S(p)$ and $m > n$, we conclude that $v > 1$, so $m \geq \overline{2}S(p)S(S(p))$. By A19, $S(S(S(p)))^2 < \overline{2}S(p)S(S(p))$ and thus $S(S(S(p)))^2 < m$. If $q$ denotes the largest prime, known to exist by A18, for which $q^2 < m$, then $q \geq S(S(S(p)))$ and the proof proceeds as in Theorem 4.1. □

## 5. On the independence of the additional axioms

That A16 and A17 are not theorems of PA$^-\cup\{$A20, A19$\}$ can be seen by noticing that $\mathcal{C}(\mathbb{Q}_\mathbb{Z}[X])$ —the positive cone of $\mathbb{Q}_\mathbb{Z}[X]$, which stands for the ring of polynomials in $X$ with constant term in $\mathbb{Z}$ and with all other coefficients in $\mathbb{Q}$, ordered by $\sum_{i=0}^n c_i X^i > 0$ if and only if $c_n > 0$ (here $c_0 \in \mathbb{Z}$, and $c_i \in \mathbb{Q}$ for all $1 \leq i \leq n$, with $c_n \neq 0$) — which is a model of PA$^-$ satisfies neither A16 nor A17, for $X + 1$, a prime, has no succesor prime and $X - 1$, also a prime, has no predecessor prime. That this is so follows from the fact that all $X + z$, with $z \in \mathbb{Z}$, with $z \neq \pm 1$, are composite.

It is plain that A19 holds in $\mathcal{C}(\mathbb{Q}_\mathbb{Z}[X])$, for it holds for numeral primes and it vacuously holds for nonstandard primes, given that there is no sequence of three primes, each of which is the successor of the previous one.

That A20 holds in $\mathcal{C}(\mathbb{Q}_\mathbb{Z}[X])$ can be seen by noticing that $f(X) = \sum_{i=0}^n c_i X^i$ in $\mathcal{C}(\mathbb{Q}_\mathbb{Z}[X])$ is a multiple of $c_0$ if $c_0 \neq \pm 1$, and $F(X)$ has as prime divisor any prime divisor of $c_0$, and that, for $c_0 = \pm 1$, $f(X)$ can be decomposed in irreducible polynomials as a polynomial in $\mathbb{Q}[X]$ with each factor having the constant term $\pm 1$.

That A18 does not follow from A17, A16, and A19 can be seen by noticing that $\mathcal{C}(\mathbb{Z}[X])$ (where all irreducibles are primes) satisfies A17, A16, and A19 but not A18. To see this, we need the following result, which is, in essence, Lemma 9 or [5, p. 234]:

**Lemma 5.1.** *A polynomial $f(X) = a_n X^n + a_{n-1} X^{n-1} + \ldots + a_1 X \pm p$, with integer coefficients and $p$ a prime number with*

$$p > \sum_{i=1}^{n} |a_i| \tag{5}$$

*is irreducible.*

*Proof.* Suppose $f(X)$ were decomposable, $f(X) = g(X)h(X)$. Both $g(X)$ and $h(X)$ would have to be polynomials of degree $\geq 1$ in $\mathbb{Z}[X]$, given that $p$ cannot divide any of the $a_i$'s. The constant term of one of $g(X)$ and $h(X)$ must be $\pm 1$, given that multiplied with the constant term of the other polynomial it becomes $\pm p$. By Viète's formulas, the product of all the complex zeros of the polynomial having constant term $\pm 1$ is $\pm \frac{1}{m}$, where $m \in \mathbb{N}$, so the absolute value of at least one of those zeros is $\leq 1$. Let $\alpha$ denote that zero. Since $\alpha$ is a zero of one of the factors of $f(X)$, $\alpha$ is a zero of $f(X)$ as well. Thus $\pm p = \sum_{i=1}^{n} a_i \alpha^i$, so $p \leq \sum_{i=1}^{n} |a_i \alpha^i| \leq \sum_{i=1}^{n} |a_i|$. This contradicts (5). □

Let $f(X) = \sum_{i=0}^{n} a_i X^i$ be an irreducible polynomial in $\mathcal{C}(\mathbb{Z}[X])$ of degree at least 1. Let $p$ be a prime number greater than $\sum_{i=0}^{n} |a_i|$. By Lemma 5.1, both $g(X) = f(X) - a_0 - p$ and $h(X) = f(X) - a_0 + p$ are irreducible polynomials in $\mathbb{Z}[X]$. Since $g(X) < f(X) < h(X)$, we conclude that A17 and A16 hold (as there must be a largest irreducible polynomial $P(f(X))$ with $g(X) \leq P(f(X)) < f(X)$ and there must be a smallest irreducible polynomial $S(f(X))$ with $f(X) < S(f(X)) \leq h(X)$). Since this obviously holds for primes that are numerals, we have shown the existence of both predecessor and successor primes.

This, incidentally, also shows that A19 holds in $\mathcal{C}(\mathbb{Z}[X])$. That it holds for primes that are numerals was shown in [2]. Let $q(X)$ be an irreducible polynomial of degree $n$ with leading coefficient $a_n$ in $\mathcal{C}(\mathbb{Z}[X])$. Then, as shown above, both $P(q(X))$ and $P(P(q(X)))$ have degree $n$ and their leading coefficients are the same as that of $q(X)$, so A19 holds as the left hand side of $q(X)^2 < \bar{2} P(q(X)) P(P(q(X)))$ is a polynomial of degree $2n$ and leading coefficient $a_n^2$ whereas the right hand side is a polynomial of degree $2n$ and leading coefficient $\bar{2} a_n^2$.

We have thus shown that

$$\text{PA}^-, \text{A19}, \text{A20} \nvdash \text{A17}$$
$$\text{PA}^-, \text{A19}, \text{A20} \nvdash \text{A16}$$
$$\text{PA}^-, \text{A16}, \text{A17}, \text{A19}, \text{A20}, \text{A21} \nvdash \text{A18}$$

To better understand the statement A19 is making, we notice that, given a prime $q \geq 3$, by A18, there exists a greatest prime $p$ such that $p^2 < \bar{2} q P(q)$. It is plain that $P(q)^2 < \bar{2} q P(q)$, as this amounts to $P(q) < \bar{2} q$. It is not as plain that the next prime, $q$ itself, is such that $q^2 < \bar{2} q P(q)$, for this is equivalent to $q < \bar{2} P(q)$, which is Chebyshev's Theorem. Moreover, A19 asks that even the next prime, $S(q)$, satisfies $S(q)^2 < \bar{2} q P(q)$, which is a stronger requirement than that of Chebyshev's Theorem, whose role, together with that of Bonse's inequality, which, again, is plainly weaker than A19, in the proof of Schatunwosky's theorem has been investigated in [4].

## 6. Final remarks

While we have seen that there are finite axiom systems, $\Sigma$ and $\Sigma'$, in which the generalized Schatunowsky theorem holds, both in its weak and its strong form, we do not know whether all the axioms that were added to $\text{PA}^-$ were actually needed.

To be more precise, let's say that an axiom $\alpha$ that is added to $\text{PA}^-$ is *needed* in the proof of a statement $\varphi$ if $\text{PA}^-, \varphi \vdash \alpha$.

First, let us show that

**Theorem 6.1.** $\mathbf{GS}^s$ *holds in* $\mathcal{C}(\mathbb{Z}[X])$.

*Proof.* We need to show only that $\mathbf{GS}^s$ holds for nonstandard primes $p$ in the hypothesis of $\mathbf{GS}^s$, for if $p$ were a numeral prime, then no nonstandard $n$ could be $p$-good, for such an $n$ would have to be divisible by all standard primes $s$ with $s > p$, given that $s^2 < n$, for otherwise $s^2$ would be a composite totative of $n$ that is co-prime with all primes $\leq p$.

The proof that $n = S(p)k_p$ is $p$-good proceeds like in Theorem 4.2, with $k_p$ defined as in that theorem, for all axioms used in that proof hold in $\mathcal{C}(\mathbb{Z}[X])$.

Suppose now $m > n$. We want to show that $m$ cannot be $p$-good. Suppose $m$ were $p$-good. Then, as in Theorem 4.2, we conclude that $m \geq 2S(p)S(S(p))$. Let $p = \sum_{i=0}^n a_i X^i$, with $a_n > 0$. Then $S(p) = b_0 + \sum_{i=1}^n a_i X^i$ and $S(S((p)) = c_0 + \sum_{i=1}^n a_i X^i$, with $a_0 < b_0 < c_0$. Thus $m$ is a polynomial of degree $d \geq 2n \geq 2$.

If $d = 2n$, then, given that, by A18, $x^2 < m$ for $x \in \{S(p), S(S(p)), S(S(S(p)))\}$, if $m$ were $p$-good, then $m$ should be divsible by $S(p), S(S(p))$, and $S(S(S(p)))$, which is not possible, as the degree of $S(p)S(S(p))S(S(S(p)))$ is $3n$. Thus no $m$ of degree $2n$ can be $p$-good.

If $d \geq 2n + 1$, then, for all $k$, $S^k(p) := S(\ldots S(p) \ldots)$, where the successor operation has been applied $k$ times, is such that $S^k(p)^2 < m$, given that the degree of the left hand side is $2n$ and that of the right hand side is $\geq 2n + 1$. We conclude that $m$ cannot be $p$-good, for it would have to be divisible by $S^k(p)$ for all $k \geq 1$. $\square$

Since $\mathcal{C}(\mathbb{Z}[X])$ satisfies all the axioms of $\Sigma'$ except A18, we have shown that

$$\text{PA}^-, \mathbf{GS}^s \nvdash \text{A18}.$$

With the meaning of an axiom being needed defined earlier, we can say that A18 is not needed in the proof of $\mathbf{GS}^s$. This does not mean that we know of a proof that does not use A18.

## References

[1] L. E. Dickson, History of the theory of numbers, vol. I (Chelsea, 1952).
[2] Y. Kaneko, H. Nakai, A generalization of Schatunowsky's theorem. Amer. Math. Monthly 132 (2025), no. 5, 443–447.
[3] R. Kaye, Models of Peano Arithmetic (Oxford University Press, 1991).
[4] V. Pambuccian, Schatunowsky's theorem, Bonse's inequality, and Chebyshev's theorem in weak fragments of Peano arithmetic. Math. Log. Quart. **61**, 230–235 (2015).
[5] H. Osada, The Galois groups of the polynomials $X^n + aX^l + b$, J. Number Theory **25**, 230–238 (1987).

School of Mathematical and Natural Studies, Arizona State University - West Valley Campus, P. O. Box 37100, Phoenix AZ 85069-7100, U.S.A.
*Email address*: hala.king@asu.edu; pamb@asu.edu